\documentclass [11pt,oneside]{amsart}

\reversemarginpar \textwidth=14cm \textheight=19cm

\usepackage{amssymb} \usepackage{amsfonts} \usepackage{amsmath}
\usepackage{amsthm} \usepackage{epsfig} 

\usepackage{caption}

\newcommand{\mettifig}[1]{\epsfig{file=#1}}

\newtheorem{lemma}{Lemma}
\newtheorem{teo}[lemma]{Theorem}
\newtheorem{prop}[lemma]{Proposition}
\newtheorem{cor}[lemma]{Corollary}

\theoremstyle{definition}
\newtheorem{defn}[lemma]{Definition}

\newtheorem{rem}[lemma]{Remark} 

\newcommand{\matN}{\ensuremath {\mathbb{N}}}
\newcommand{\matR} {\ensuremath {\mathbb{R}}}

\newcommand{\matZ} {\ensuremath {\mathbb{Z}}}
\newcommand{\matC} {\ensuremath {\mathbb{C}}}
\newcommand{\matP} {\ensuremath {\mathbb{P}}}
\newcommand{\matH} {\ensuremath {\mathbb{H}}}

\newcommand{\Isom}{\ensuremath{{\rm Isom}}}
\newcommand{\MCG}{\ensuremath{{\mathcal M \mathcal C \mathcal G}}}
\newcommand{\Out}{\ensuremath{{\rm Out}}}
\newcommand{\Aut}{\ensuremath{{\rm Aut}}}

\def\freccinabis#1#2{\begin{picture}(40,15)
\put(0,-5){\mettifig{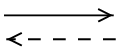}}
\put(#2,-13){\tiny #1}\end{picture}}

\def\freccinavertbis#1{\begin{picture}(20,20)
\put(0,0){\mettifig{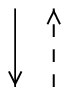}}
\put(28,10){\tiny #1}\end{picture}}

\author{Roberto Frigerio}
\author{Bruno Martelli}

\address{Dipartimento di Matematica \\%
Universit\`a di Pisa \\%
Largo B.~Pontecorvo 5 \\%
56127 Pisa, Italy%
}

\email{frigerio@mail.dm.unipi.it, martelli@dm.unipi.it}

\title[Countable groups are mapping class groups]{Countable groups are
mapping class groups\\ of hyperbolic $3$-manifolds}

\subjclass[2000]{57M50 (primary), 20F29, 30F40 (secondary).}

\keywords{Kleinian groups, hyperbolic $3$-manifolds, special polyhedra.}

\thanks{Both authors are partially supported by the INTAS project
``CalcoMet-GT'' 03-51-3663}

\begin{document}

\begin{abstract}
We prove that for every 
countable group $G$ there exists a hyperbolic $3$-manifold $M$ 
such that the isometry group of $M$, the mapping class group of $M$, and the outer 
automorphism group of $\pi_1 (M)$ are isomorphic to $G$.
\end{abstract}

\maketitle

\section*{Introduction}

A \emph{hyperbolic manifold} here is a connected orientable
paracompact manifold (without boundary) equipped
with a complete metric of constant sectional curvature equal to $-1$. 
If $M$ is a hyperbolic $3$-manifold,  $\Isom (M)$ is
the group of isometries of $M$, and $\MCG(M)$
the \emph{mapping class group} of $M$, \emph{i.e.}~the 
group of isotopy classes of  
self-homeomorphisms of $M$. 
If $G$ is a group, ${\rm Out} (G)$ is the outer automorphism
group of $G$. We prove here the following.

\begin{teo}\label{main2:teo}
For every countable group $G$ there is a hyperbolic $3$-manifold $M$
such that: 
\[
G\cong\Isom(M)\cong \MCG (M)\cong\Out (\pi_1 (M)).
\]
\end{teo}

Together with the fact that there exist uncountably many 
pairwise non-isomorphic countable groups, Theorem~\ref{main2:teo}
implies the following:

\begin{cor}\label{main1:cor}
There are uncountably many non-isomorphic groups that are fundamental groups
of hyperbolic $3$-manifolds.
\end{cor}

Every hyperbolic $3$-manifold
is the quotient of $\matH^3$ via the action of a \emph{Kleinian group}, that is
a discrete torsion-free subgroup of $\matP {\rm SL}_2(\matC)$.
Theorem~\ref{main2:teo} implies the following:

\begin{cor}\label{group:cor}
Every countable group is the outer automorphism group of a
Kleinian group.
\end{cor}

\subsection*{Related results}
Theorem~\ref{main2:teo} is already known for finite groups:
Kojima proved in~\cite{Kojima} that every finite
group is the isometry group of a compact hyperbolic $3$-manifold
(see also \cite{CoFriMaPe}).
Moreover, if $M$ is compact hyperbolic, 
Mostow's rigidity Theorem \cite{Mos} and
Gabai-Meyerhoff-Thurston's 
%homotopy-implies-isotopy
%Theorem 
results~\cite{Gabai}
imply that $\Out (\pi_1 (M))$, $\Isom (M)$, 
and $\MCG (M)$ are isomorphic finite groups.
Kojima's Theorem for finite groups has been
extended recently to higher dimensions~\cite{Kojima+}, while
Winckelmann proved in~\cite{surfaces} that every countable group
is the isometry group of a complete hyperbolic surface.

There are many analogues of Corollary~\ref{main1:cor} concerning 
topological manifolds
of various dimensions.
In dimension $2$, there are uncountably many (pairwise non-homeomorphic) 
paracompact surfaces~\cite{giapponese}. 
On the other hand
every paracompact surface
supports a complex structure, so Riemann's uniformization Theorem 
and classical results on Fuchsian 
groups~\cite{Fuchs2,Fuchs} imply that there exist only
countably many isomorphism classes of fundamental groups
of surfaces.   
In dimensions higher than $2$,
there are uncountably many contractible 
manifolds~\cite{contractible,contractible4,contractible5}, 
and also uncountably many fundamental groups
(see \emph{e.g.}~\cite{contractible4,groups}).

Regarding Corollary~\ref{group:cor},
Matumoto proved in~\cite{Matumoto} that every group 
is the outer automorphism group of some group, 
while Bumagin and Wise recently showed that every
countable group is the outer automorphism group
of a two generators group~\cite{groupout}.  

\subsection*{Sketch of the proof}
We consider some decorated $2$-dimensional polyhedra with infinitely many vertices,
edges, and faces. We associate to each such polyhedron a hyperbolic $3$-manifold
by assembling some blocks with geodesic boundary, corresponding
to vertices and edges of the polyhedron. We show that the isometry
group of the manifold is isomorphic to the combinatorial automorphism group
of the polyhedron. It is easy to see that every countable group is
the combinatorial automorphism group of one such polyhedron.

If $M$ is a hyperbolic manifold arising from our construction, building on
a rigidity result of Keen, Maskit, and Series~\cite{KMS} 
we show that every automorphism of $\pi_1(M)$ is induced
by an isometry. This also implies that a self-homeomorphism of $M$ is homotopic to an isometry:
it remains to show that it is indeed isotopic to it. This follows from a
result of Brown~\cite{brown2}, since our $M$ is
non-compact and \emph{end-irreducible}.

The paper is organized as follows. We describe the decorated polyhedra
and their associated hyperbolic manifolds in Section~\ref{construction:section}.
A more detailed outline of the proof of Theorem~\ref{main2:teo} is then given
in Section~\ref{scheme:section}.
Polyhedra with assigned automorphism group are constructed in
Section \ref{special:section}. The needed rigidity results are then
proved in Section~\ref{maps:section}.

\newpage

\section{From decorated polyhedra to hyperbolic manifolds}\label{construction:section}

\subsection{Decorated polyhedra}\label{construction:subsection}
\begin{figure}
\begin{center}
\mettifig{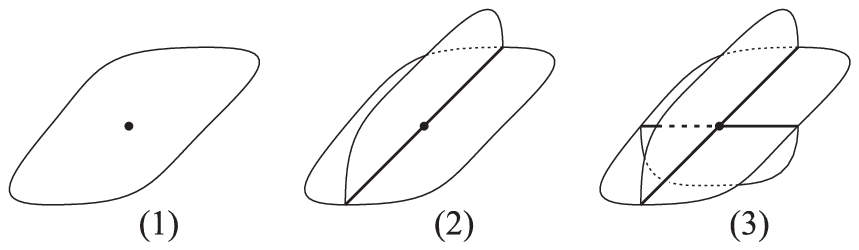} 
\caption{Neighbourhoods of points in a special polyhedron.}
\label{new_standard_nhbds:fig}
\end{center}
\end{figure}

A $2$-dimensional polyhedron $P$ is \emph{special} if every point
of $P$ has one of the regular neighbourhoods shown in 
Fig.~\ref{new_standard_nhbds:fig}, and if 
the stratification given by the three
types of points gives a cellularization of $P$. That is, points
of type (1) form discs -- the \emph{faces} -- and points of type (2) form segments
-- the \emph{edges}.
A \emph{vertex} is then a point of type (3).
In contrast with the usual definition, compactness is not required here:
thus $P$ can have infinitely many vertices, edges, and discs
(but each has compact closure and is incident to finitely many other
vertices, edges, and discs).

A \emph{decorated polyhedron} is a pair $(P,c)$,
where $P$ is a special polyhedron and $c:\{ {\rm edges\ of}
\ P\}\to \matN$ is a function, which will be called the
\emph{decoration} of $P$. 
The aim of this section is to describe a recipe for
associating to every decorated polyhedron $(P,c)$ 
an orientable hyperbolic $3$-manifold $M(P,c)$. This will be done by
suitably modifying
the \emph{shadow} construction first introduced 
by Turaev and subsequently explored in~\cite{CoTh, CoFriMaPe}.

\subsection{From polyhedra to hyperbolic manifolds}
\begin{figure}
\begin{center}
\mettifig{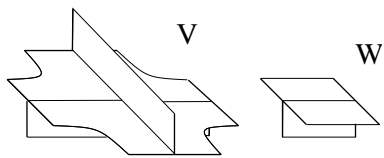, width = 5 cm} 
\caption{The pieces $V$ and $W$.}
\label{pezzi:fig}
\end{center}
\end{figure}

Let $(P,c)$ be a decorated polyhedron.
Let $P_0$ be a regular neighborhood of its $1$-skeleton.
We subdivide $P_0$ into copies of the pieces shown in Fig.~\ref{pezzi:fig}: we take one
piece of type $V$ for each vertex and $c(e)$ pieces of type $W$ for
each edge $e$. 
The ``boundary'' of $V$ (resp.~$W$) consists of $4$ (resp.~$2$) $Y$-shaped
graphs and $6$ (resp.~$3$) arcs. The following result will be proved
in Lemmas~\ref{isoK:lemma} and~\ref{E:lemma} below.

\begin{prop}\label{grosso:prop}
There are two complete finite-volume hyperbolic $3$-manifolds 
$K$ and $E$ with geodesic boundary such that:
\begin{itemize}
\item
their boundary consists of thrice-punctured spheres;
\item
boundary components of $K$ (resp.~of $E$) 
correspond to $Y$-shaped graphs in $\partial V$ (resp.~in $\partial W$),
with punctures corresponding to endpoints of the graphs;
\item
annular cusps of $K$ (resp.~of $E$) 
correspond to arcs in $\partial V$ (resp.~in $\partial W$);
\item
the intersections between 
the boundary components and the annular cusps of $K$ (resp.~of $E$) 
correspond to the intersections between the corresponding
$Y$-shaped graphs and arcs in $\partial V$ (resp.~in $\partial W$);
\item
the actions of the isometries on boundary components and cusps induce the isomorphisms:
\[
\begin{array}{c}
\Isom^+(K) \cong \Aut(\partial V) \cong \mathfrak{S}_4,\\
\Isom(E) \cong \Isom^+(E) \cong \Aut(\partial W) \cong \mathfrak{S}_3 \times \matZ/_2,
\end{array}
\]
\item
every geodesic thrice-punctured sphere in $K$ (resp.~in $E$) is contained in $\partial K$
(resp.~in $\partial E$).
\end{itemize}
\end{prop}

We have denoted by $\Aut(Z)$ the combinatorial automorphism group of a graph $Z$, and
by $\Isom^+$ the group of orientation-preserving isometries.

We equip the blocks $K$ and $E$ with arbitrary orientations.
We construct from $(P,c)$ an oriented hyperbolic $3$-manifold
as follows: as we said above, the polyhedron $P_0$ decomposes 
into pieces of type $V$ and $W$. 
Pick blocks of type $K$ and $E$ corresponding to pieces of type $V$ and $W$.
For each $Y$-shaped graph we have a thrice-punctured sphere, whose
punctures correspond to the endpoints of the graph.
Given two geodesic thrice-punctured spheres, 
every bijection between their punctures
is realized by a unique \emph{orientation-reversing} isometry.
The identifications of the $Y$-shaped graphs in $P_0$ therefore induce
isometries between pairs of boundary punctured spheres, and
we use such isometries for gluing all these pairs.

Since any symmetry of $V$ (resp.~of $W$) translates into
an orientation-preserving isometry of $K$ (resp.~of $E$),
there is no ambiguity in the construction, and
the result is a complete oriented hyperbolic $3$-manifold $M(P,c)$.
Annular cusps glue up to toric cusps, which are in natural correspondence 
with the faces of $P$.
Every vertex of $P$ corresponds to a $K$-block.
Every edge $e$ of $P$ gives rise 
to $c(e)$ blocks of type $E$ and $c(e)+1$ thrice-punctured spheres.

\begin{rem}\label{exh:rem}
If $(P,c)$ is a decorated 
polyhedron, there are connected complete hyperbolic $3$-submanifolds
with geodesic boundary $M_i\subset M(P,c)$, $i\in\matN$ with
the following properties:
\begin{itemize}
\item
$M_i$ is the union of a finite number
of $K$-blocks and $E$-blocks,
for all $i\in\matN$;
\item
$M_i\subset {\rm int}\, M_{i+1}$ for all $i\in\matN$, and
$\bigcup_{i=0}^\infty M_i=M(P,c)$.
\end{itemize}
\end{rem}

\subsection{The block $K$}
The rest of this section is devoted to the proof of Proposition~\ref{grosso:prop}. 
The block $K$ we introduce here has already been used by various authors
to construct hyperbolic manifolds. As we mentioned above, Turaev has studied
the boundary of the $4$-dimensional thickening of a special polyhedron
(a \emph{shadow}), and subsequently Costantino and Thurston have decomposed it
into $K$-blocks and solid tori, lying above 
the vertices and the faces of the polyhedron~\cite{CoTh}.
The block $K$ has also been used by
Minsky to construct the \emph{model manifold}
needed in the proof of the Ending Lamination Conjecture~\cite{Minsky}, and
by Agol to construct non-Haken manifolds of arbitrarily high genus~\cite{Agol03}.

\begin{figure}
\begin{center}
\mettifig{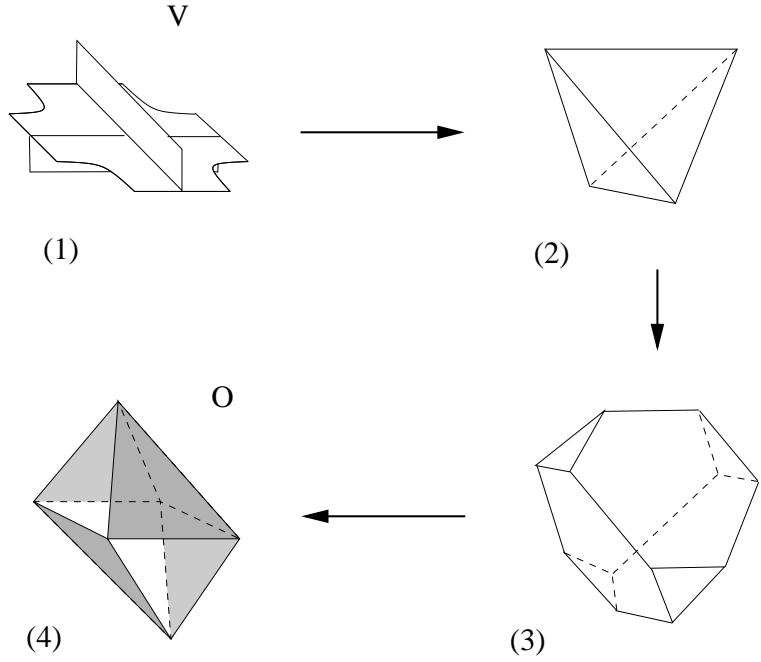, width = 10cm} 
\caption{The polyhedron $P_0$ is made of pieces whose boundary consists
of $4$ $Y$-shaped graphs and $6$ arcs (1).
Each such piece is combinatorially equivalent to a truncated tetrahedron
(3), whence to a regular ideal octahedron with checkerboard coloring
of the faces (4).}
\label{crea7:fig}
\end{center}
\end{figure}

The piece $V$
is dual to a tetrahedron, which is in turn combinatorially
equivalent to an ideal octahedron $O$ with a checkerboard coloring of the faces:
see Fig.~\ref{crea7:fig}. Note that under this correspondence,
the $Y$-shaped graphs and the six arcs in $\partial V$
are associated respectively to the shadowed faces and the 
vertices of the octahedron.

\begin{figure}
\begin{center}
\mettifig{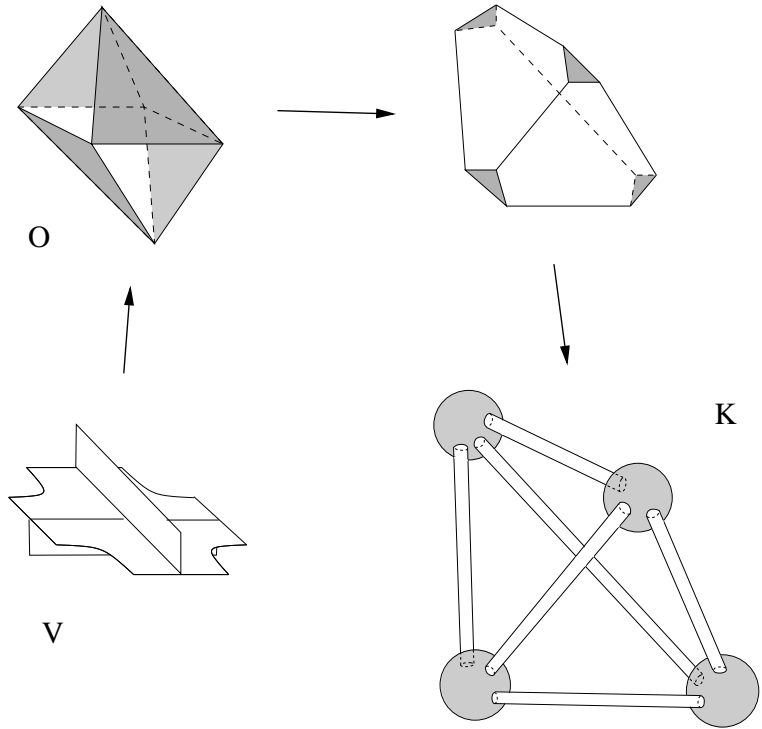, width = 10cm} 
\caption{Each piece $V$ gives rise to a block $K$ which is homeomorphic
to the complement
in $S^3$ of four open balls connected by six arcs. The thrice-punctured spheres
in $\partial K$ correspond to the $Y$-shaped graphs in $\partial V$.}
\label{crea6:fig}
\end{center}
\end{figure}

Now let us realize $O$ as an ideal regular hyperbolic 
octahedron, and define $K$ as the geometric object obtained by mirroring
$O$ along its white faces. Since the dihedral
angles of the regular ideal octahedron are right, $K$ is 
a complete hyperbolic manifold with non-compact geodesic boundary. 
As required, $\partial K$ consists of four geodesic thrice-punctured 
spheres, each of which canonically corresponds to a $Y$-shaped graph
in $\partial V$. The six vertices of $O$ give rise to six annular cusps
of $K$, which in turn intersect each component of $\partial K$ in
(neighbourhoods of) the punctures. Thus a puncture on a component
of $\partial K$ canonically corresponds to the intersection point of 
one $Y$-shaped graph with an arc on $\partial V$ (see Fig.~\ref{crea6:fig}).

\begin{lemma}\label{isoK:lemma}
We have $\Isom (K)\cong\mathfrak{S}_4\times \matZ/_2$. The isomorphism
is obtained by associating to each isometry the induced permutation of the components
of $\partial K$ and a number in $\{0,1\}$ saying whether 
the isometry is orientation-preserving or not.
The isometry corresponding to $(id, 1)$
switches the ideal octahedra in $K$. 

If $S\subset K$ is an embedded geodesic thrice-punctured sphere,
then $S$ is a boundary component of $K$.
\end{lemma}
\begin{proof}
The map from $\Isom(K)$ to $\mathfrak{S}_4\times \matZ/_2$ is surjective, because
$K$ inherits the symmetries of the tetrahedron. Concerning injectivity, an
isometry in the kernel fixes each thrice-punctured sphere in $\partial K$, 
each of its punctures, and is orientation-preserving: therefore it fixes
$\partial K$ pointwise, and hence the whole of $K$.

Miyamoto proved in~\cite{Miy} that if $N$ is a complete finite-volume
hyperbolic manifold with geodesic boundary, then
$2{\rm vol} (N)\geqslant -\chi (\partial N)\cdot v_O$, where
$v_O$ is the volume of the regular ideal octahedron.
Suppose that $S\subset K\setminus \partial K$ 
is an embedded thrice-punctured sphere, and let $N$ be the 
(possibly disconnected) hyperbolic manifold with geodesic boundary
obtained by cutting $K$ along $S$.  We have 
$4v_O=2{\rm vol} (K)=2{\rm vol} (N)\geqslant -\chi (\partial N) v_O=
6v_O$, a contradiction. 
\end{proof}

\subsection{The block $E$}
Let $E_0$ be  
the manifold obtained by cutting the double of $K$ along
one component of $\partial K$ (by Lemma~\ref{isoK:lemma},
the isometry type of this manifold does not depend
on the choice of the component of $\partial K$ along which 
we cut).
By construction $E_0$ decomposes
into two $K$-blocks, has $3$ toric cusps, and $3$ annular cusps.
Moreover, $\partial E_0$ is given by two 
geodesic thrice-punctured spheres.

\begin{lemma}\label{E0:lemma}
We have $\Isom(E_0)\cong\mathfrak{S}_3\times\matZ/_2\times\matZ/_2$.
The isomorphism is obtained by associating to each isometry the induced permutation
on the three annular (or toric) cusps, on its two boundary components, 
and a number in $\{0,1\}$ saying whether it is orientation-preserving or not.

Moreover, if $S$ is a geodesic thrice-punctured sphere embedded in
$E_0\setminus \partial E_0$, then at least one puncture of $S$
lies in a toric cusp of $E_0$.
\end{lemma}
\begin{proof}
The map from $\Isom(E_0)$ to $\mathfrak{S}_3\times\matZ/_2\times\matZ/_2$
is surjective by Lemma~\ref{isoK:lemma}, and is injective because
an isometry in the kernel fixes $\partial E_0$ pointwise.

Suppose $S\subset E_0\setminus\partial E_0$ is an embedded
geodesic thrice-punctured sphere not intersecting the
toric cusps of $E_0$. Then $S$ should intersect the punctured spheres
separating the two blocks in closed geodesics. 
Since there are no closed geodesics in thrice-punctured spheres,
this would imply that $S$ is contained 
in the interior of one of the two $K$-blocks,
against Lemma~\ref{isoK:lemma}.
\end{proof}

The block $E_0$ fullfills all the requirements to be our block $E$,
except two: it has an orientation-reversing involution
and contains some thrice-punctured spheres (the ones separating the two blocks).
We now kill all the redundant isometries and punctured spheres
by filling appropriately the three toric cusps.

We recall that a \emph{slope} on a torus is an isotopy class
of simple closed unoriented curves. 
Take three disjoint sections of the toric cusps of the same area:
their boundaries are Euclidean tori, whose slopes have a definite length.
The group ${\rm Isom}^+ (E_0)$ acts on the set of all slopes on the $3$
toric cusps. A direct analysis shows that 
the orbit of any slope $s$ consists of $3$ slopes of the same length, 
each lying on a distinct toric cusp. 

An easy doubling argument~\cite{FujKoj} shows that
Thurston's hyperbolic Dehn filling Theorem 
also applies to hyperbolic manifolds with geodesic boundary.
Take then a sufficiently long slope $s$ on a toric cusp, such that the following 
conditions hold:
\begin{enumerate}
\item
the slope $s$ is neither ``vertical'' or ``horizontal'', so the
corresponding triple is not invariant under orientation-reversing isometries;
\item
on the chosen horospherical sections of the toric cusps,
the slopes  corresponding to 
$s$ have length $\ell(s)> 12$;
\item
the manifold $E$, obtained via Dehn filling the $3$ toric cusps by 
killing the $3$ slopes corresponding to $s$,
is hyperbolic with geodesic boundary;
\item
the shortest closed geodesics in $E$ are the cores
of the added solid tori.
\end{enumerate}

\begin{lemma}\label{E:lemma}
Every positive isometry of $E_0$ induces an isometry of $E$, and this gives the isomorphism
${\rm Isom} (E)=\Isom^+ (E)\cong \mathfrak{S}_3\times\matZ/_2$.
Moreover, the only geodesic thrice-punctured spheres contained in $E$ are the components
of $\partial E$.	
\end{lemma}
\begin{proof}
Every positive isometry of $E_0$ leaves the triple of killed slopes invariant, and hence
extends to a self-homeomorphism of $E$, which is homotopic to an isometry of $E$ by Mostow
rigidity for manifolds with boundary~\cite{FriPe}.
Conversely, every isometry of $E$ fixes the set of the $3$ shortest slopes, and hence
its complement $E_0$. Therefore it is homotopic to an isometry of $E_0$ fixing
the triple of slopes, which must be orientation-preserving by assumption (1) above.

Suppose we have a geodesic thrice-punctured sphere $S$ in $E$. 
Since $S$ does not contain closed geodesics, it
intersects the three added geodesics in some $k$ points, so
$S_0=S\cap E_0$ is a sphere with $3+k$ punctures inside $E_0$.
Since $S$ is geodesic, it is easily seen that
$S_0$ is \emph{essential}, that is
it is incompressible and $\partial$-incompressible, in the sense of \cite{Agol}.
Therefore \cite[Theorem 5.1]{Agol}
(or equivalently \cite{Lack}) implies that
$12k<k\ell(s)\leqslant |6\chi (S)|=6(1+k)$, whence $k=0$.
Thus $S=S_0$ lies in $E_0$. By~\cite{Adams}, $S$ is isotopic
to an embedded geodesic thrice-punctured sphere in $E_0$ not intersecting
the toric cusps of $E_0$, whence,
by Lemma~\ref{E0:lemma},
to a component of $\partial E_0$.
This readily implies that $S$ is a component of $\partial E$.
\end{proof}

\begin{rem}\label{orientation:rem}
If only the block $K$ were used (as in~\cite{CoFriMaPe}), then
the orientation-reversing
involution switching the ideal octahedra of this block
would always extend to the whole of $M(P,c)$, thus
giving an annoying $\matZ/_2$ factor in $\Isom (M(P,c))$.
In~\cite{CoFriMaPe}, this factor was killed via Dehn filling, but
Thurston's hyperbolic Dehn filling Theorem cannot be applied
easily  here because the manifold is very big.
This is the main reason for introducing
the block $E$, which is 
\emph{chiral}, \emph{i.e.}~it has no orientation-reversing
isometries.
The introduction of $E$ also provides
some technical semplifications when constructing a polyhedron
with the appropriate combinatorial isomorphism group.
\end{rem}

\section{Scheme of the proof of Theorem~\ref{main2:teo}} \label{scheme:section}
We say that a special polyhedron is \emph{regular}
if it has more than two vertices, and 
every edge in its $1$-skeleton has distinct endpoints.
A decorated polyhedron $(P,c)$ is \emph{big} if 
$P$ is regular and $c(e)\geqslant 4$ for every
edge $e$ of $P$.
We denote by ${\rm Aut} (P,c)$ the group of
combinatorial automorphisms of $(P,c)$, \emph{i.e.}~the group
of  combinatorial automorphisms of $P$ preserving the decoration.

The proof of Theorem~\ref{main2:teo} consists of two steps: 
\begin{enumerate}
\item
we show that
every countable group is the group of combinatorial automorphisms of
some big decorated polyhedron: 
this is done in Section~\ref{special:section};
\item
for any big decorated polyedron $(P,c)$,
we prove in Section~\ref{maps:section} 
bijections of all the following sets:
\end{enumerate}

%\begin{center}
%\begin{tabular}{@{}c@{}c@{}c@{}c@{}c@{}} 

%$\left\{\!\begin{array}{c} {\rm combinatorial} \\ {\rm automorphisms} \\ 
%{\rm of\ } (P,c) \end{array}\!\right\}$ &

%\ \ \freccina{(Prop.~\ref{isom2comb:prop})}2\!\!\!\! &

%$\left\{\!\begin{array}{c} {\rm isometries} \\ {\rm of\ } M(P,c)
%\end{array}\!\right\}$ &

%\freccina{(Th.~\ref{group:isometry:teo})}4 &

%$\left\{\!\begin{array}{c} {\rm outer} \\ {\rm automorphisms} \\ 
%{\rm of\ }\pi_1(M(P,c))\end{array}\!\right\}$ \\

% & & \freccinavert{(Cor.~\ref{cor:cor})} & & \\

%$\left\{\!\begin{array}{c} \textrm{self-homeom.} \\ {\rm of\ } M(P,c)\\ {\rm\ up\ to\  isotopy}  
%\end{array}\!\right\}$ &

%\!\freccina{(Prop.~\ref{homotopy2isotopy})}0 &

%$\left\{\!\begin{array}{c} \textrm{self-homeom.} \\ {\rm of\ } M(P,c)\\
%{\rm\ up\ to\ homotopy} 
%\end{array}\!\right\}$ &

%\end{tabular}
%\end{center}

\begin{center}
\begin{tabular}{@{}c@{}c@{}c@{}c@{}c@{}} 

${\rm Aut} (P,c)$ &

 \freccinabis{(Prop.~\ref{isom2comb:prop})}2 &

$\Isom (M(P,c))$ &

\freccinabis{(Th.~\ref{group:isometry:teo})}4 &

\ ${\rm Out} (\pi_1 (M(P,c)))$ \\

& & & & \\
 & & \freccinavertbis{(Cor.~\ref{cor:cor})} \ \ \ \ \ & & \\

$\MCG (M(P,c))$\ \ \  &

\freccinabis{(Prop.~\ref{homotopy2isotopy})}0 &

$\MCG^{\rm hom} (M(P,c))$  &

\end{tabular}
\end{center}

\vspace {.5 cm}

We have denoted by $\MCG^{\rm hom} (M(P,c))$ the group of 
homotopy classes of self-homeomorphisms
of $M(P,c)$. 
Maps corresponding to solid arrows are natural, while their inverses,
corresponding to dashed arrows, are constructed in Section~\ref{maps:section}.

\section{Decorated polyhedra with assigned automorphism group}\label{special:section}
We prove here that every countable group is the automorphism group of a big
decorated polyhedron.

\subsection{Special polyhedra with assigned fundamental groups}
We begin with the following:
\begin{prop}\label{construct:P:prop}
Every countable group is the fundamental group of a regular special polyhedron.
\end{prop}
\begin{proof}
A countable group has some presentation with countably many generators and relators.
In general, the corresponding
polyhedra are constructed by attaching a $1$-cell for each generator
and a $2$-cell for each relator to some simply connected base-space (for instance, a point). 
Here, since cells need to be locally finite, the base-space needs to be non-compact.

Our base-space is the $2$-skeleton of an infinite wall, 
\emph{i.e.}~of a
strip $\matR^2\times [0,1]\subset\matR^3$  tessellated 
as in Fig.~\ref{bricks:fig}. 
It is a simply connected special polyhedron
whose horizontal faces are hexagons
and vertical faces are squares.
\begin{figure}
\begin{center}
\mettifig{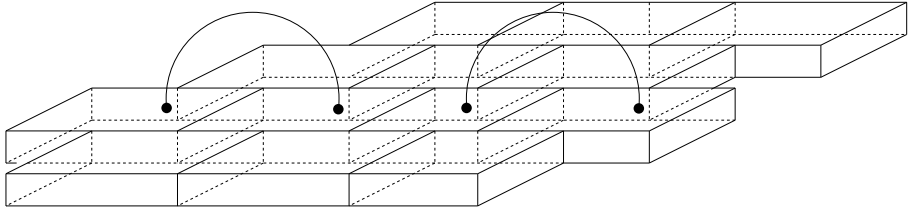, width = 12cm} 
\caption{The $2$-skeleton of an infinite wall is a special polyhedron, and we add
$1$-cells to pairs of consecutive bricks in an infinite line.}
\label{bricks:fig}
\end{center}
\end{figure}
Take an infinite line of consecutive bricks, a point inside the top-face of each
brick, and add $1$-cells to pairs of consecutive points, as in Fig.~\ref{bricks:fig}. 

In our presentation, 
we can suppose that every generator occurs in exactly $3$ relators, 
and once in each~\cite{relatori}\footnote{For each generator $g$, 
substitute its occurrences in the relators with new distinct generators
$g_1,\ldots,g_k$, each occurring once, and add new relators $g_ig_{i+1}^{-1}$.}.
Assign arbitrarily an orientation and a generator to each
$1$-cell, and attach for each relator 
a $2$-cell running alternatively along a $1$-cell
(corresponding to a letter) and as a right-angled arc along the wall as in Fig.~\ref{bricks2:fig}.
\begin{figure}
\begin{center}
\mettifig{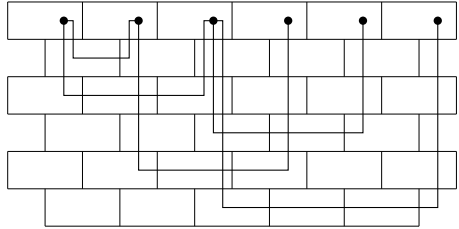, width = 6cm} 
\caption{The $2$-cells run over the bricks along right-angled paths.}
\label{bricks2:fig}
\end{center}
\end{figure}
Each new two-dimensional region is a cell, and
three distinct $2$-cells run along each $1$-cell. 
Therefore the resulting polyhedron is special. It is easy to see that it is regular.
\end{proof}

\begin{prop}
Every countable group is the group of combinatorial
automorphisms of a big decorated polyhedron.
\end{prop}
\begin{proof}
In a decorated polyhedron, let us define the colour of a vertex as the $4$-uple
of colours of the adjacent edges.
Let $G$ be a countable group and $P$ be a regular special polyhedron 
with fundamental group $G$.

Let $e_0,e_1,\ldots,e_n,\ldots$ be an arbitrary ordering of the edges of $P$,
and set $c(e_i)=i+4$ for all $i\in\matN$. 
Since $P$ is regular, distinct vertices have distinct colours
(because edges have distinct endpoints, and the singular locus of $P$
is not \mettifig{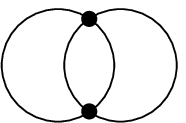, width=.5cm}).
 
Let now $\widetilde{P}$ be the universal cover of $P$, and let 
$\widetilde{c}$ be the decoration of $\widetilde{P}$ induced by $c$.
The fiber over a vertex or edge consists precisely of all vertices or edges
of the same colour.
Of course the group of
deck transformations of $\widetilde P$ is isomorphic to $G$.
By definition of $\widetilde{c}$, deck transformations are automorphisms.
Conversely, every automorphism is a deck transformation, because
it must preserve the colours of vertices and edges, and hence the fibers.
\end{proof}

\section{Rigidity results}\label{maps:section}
\subsection{From isometries to combinatorial homeomorphisms of polyhedra}
We begin with the following:
\begin{prop}\label{spheres:prop}
Let $(P,c)$ be a big decorated polyhedron.
Then every geodesic thrice-punctured
sphere in $M(P,c)$ is a boundary component of some block.
\end{prop}
\begin{proof}
Suppose a geodesic thrice-punctured sphere $S$ is not a boundary component
of a block. By Proposition~\ref{grosso:prop}, it cannot be contained inside a block,
and hence it intersects some boundary component $S'$ of some block in a geodesic.
The only simple geodesics 
available in $S$ are lines connecting punctures. Thus $S$ intersects
$S'$ along a cusp, and it must intersect also every other punctured sphere which
is parallel to $S'$ along that cusp. Since $(P,c)$ is big, there are
at least $4$ such punctured spheres: a contradiction,
since there are at most $3$ pairwise disjoint geodesics on $S$.
\end{proof}

\begin{prop}\label{isom2comb:prop}
Let $(P,c)$ be a big decorated polyhedron. Then
$$\Isom(M(P,c))=\Isom^+ (M(P,c))\cong {\rm Aut} (P,c).$$
\end{prop}
\begin{proof}
Every isometry of $M(P,c)$ preserves the set of all geodesic 
thrice-punctured spheres, whence, by Proposition~\ref{spheres:prop},
the decomposition of $M(P,c)$ into
$K$-blocks corresponding to vertices and $E$-blocks corresponding to edges.
Together with
Proposition~\ref{grosso:prop}, this implies that
the group of combinatorial
automorphisms of $(P,c)$ is canonically isomorphic to the group
of \emph{orientation-preserving} isometries
of $M(P,c)$, 
so we are left to prove that $\Isom (M(P,c))=\Isom^+ (M(P,c))$.

Let  $\varphi\in\Isom(M(P,c))$
and $E_1,E_2\subset M(P,c)$
be $E$-blocks with $\varphi (E_1)=E_2$. By construction, for $i=1,2$
there exists an
orientation-preserving isometry $\psi_i$ between $E_i$ and the standard
block $E$. If $\varphi$ were orientation-reversing, then 
$\psi_2\varphi\psi_1^{-1}$ would provide an orientation-reversing
element of $\Isom (E)$, against the chirality of $E$.
\end{proof}

\subsection{Hyperbolic manifolds with 
geodesic boundary}\label{hyperbolic:subsection}
Let $N$ be a complete finite-volume hyperbolic $3$-manifold
with non-empty geodesic boundary. Its universal covering
$\widetilde{N}$ 
is isometric to
a convex polyhedron of $\matH^3$ bounded by
a countable number of disjoint geodesic hyperplanes~\cite{Kojima1}.
The group of covering automorphisms of
$\widetilde{N}$ is 
a Kleinian group $\Gamma$ acting on
$\widetilde{N}$ with
$N\cong\widetilde{N}/\Gamma$.
Kojima showed in~\cite{Kojima1} that $N$ is the \emph{convex core}
of $\Gamma$ (see~\cite{Thu} for a definition).

\subsection{Maximally parabolic Kleinian groups}
Let $G$ be a finitely generated Klein\-ian group and set 
$M=\matH^3/G$. Let $M_c\subset M$ be a compact core of $M$~\cite{Scott}
and set $b(G)=-3\chi (M_c)$.
It is shown in~\cite{KMS} that $b(G)$ only depends on the isomorphism
type of $G$ as an abstract group, and that
$G$ contains at most $b(G)$ conjugacy classes of rank-$1$
maximal parabolic subgroups.
A finitely generated Kleinian group $G$ containing $b(G)$ conjugacy classes of
rank-$1$ maximal parabolic subgroups is called
\emph{maximally parabolic}.
The following result is a restatement of~\cite[Theorems I and II]{KMS},
and gives a complete
characterization of maximally parabolic Kleinian groups.
\begin{teo}\label{maskit1:teo}
A finitely generated Kleinian group $G$ 
is maximally parabolic if and only if its
convex core is a complete finite-volume hyperbolic
$3$-manifold with geodesic boundary, consisting of a finite number of
totally geodesic thrice-punctured spheres. 
\end{teo}
Let $G$ and $G'$ be Kleinian groups and suppose 
$\varphi:G\to G'$ is an isomorphism. We say that
$\varphi$ is \emph{type-preserving} if
$\varphi$ sends parabolic elements to parabolic elements
and loxodromic elements to loxodromic elements.
The following rigidity theorem for maximally parabolic Kleinian groups
is taken from~\cite{KMS}.
\begin{teo}\label{maskit2:teo}
Let $\varphi:G\to G'$ be a type-preserving isomorphism between two
Kleinian groups $G$ and $G'$. If $G$ is maximally parabolic,
then there exists an element $h\in{\rm Isom}
(\matH^3)$ such that $\varphi (g)=h g h^{-1}$
for all $g\in G$.
\end{teo}

\subsection{From isomorphisms of fundamental groups to isometries}
Let $(P,c)$ be any decorated polyhedron, and let
$\Gamma$ 
be the Klenian group with
$M(P,c)\cong\matH^3/\Gamma$. 

\begin{teo}\label{group:isometry:teo}
Let $\varphi:\Gamma\to\Gamma$ be an isomorphism of abstract groups.
Then 
there exists an isometry $g\in\Isom (\matH^3)$ with
$\varphi (\gamma)=g\gamma g^{-1}$ for all $\gamma\in\Gamma$.
\end{teo}
\begin{proof}
If $P$ is compact then
the conclusion
follows from Mostow-Prasad's
rigidity Theorem, so
we concentrate here on the case when $P$ is 
non-compact.

The parabolic elements of $\Gamma$
can be characterized as those elements belonging to a $\matZ\oplus\matZ$
subgroup of $\Gamma$, so
$\varphi$ is type-preserving.

Let $M_i\subset M(P,c)$, $i\in\matN$ be the finite-volume
manifolds with geodesic boundary described in Remark~\ref{exh:rem},
and take a basepoint $x_0\in M_0$.
The map $\pi_1 (M_i,x_0)\to\pi_1 (M(P,c),x_0)$ 
induced by the inclusion
is injective, because $M_i$ has geodesic (and hence incompressible) boundary.
Let $\Gamma_i<\Gamma$
be the subgroup corresponding to $\pi_1 (M_i,x_0)$ under the identification
$\pi_1 (M(P,c),x_0)\cong \Gamma$.
By Theorem~\ref{maskit1:teo},
$\Gamma_i$ is maximally parabolic, so
Theorem~\ref{maskit2:teo} implies that
for every $i\in\matN$ there exists $g_i\in{\rm Isom} (\matH^3)$
such that $\varphi (\gamma)=g_i\gamma g_i^{-1}$
for all $\gamma\in\Gamma_i$. 
It follows that for all $i\in\matN$ the isometry
$g_i g_0^{-1}$ commutes with all the elements in $\Gamma_0$.
Since $\Gamma_0$ is non-elementary, this implies
$g_0=g_i$ for all $i\in\matN$, whence
$\varphi (\gamma)=g_0 \gamma g_0^{-1}$
for every $\gamma\in\bigcup_{i\in\matN} \Gamma_i=\Gamma$.
\end{proof}

\begin{cor} \label{cor:cor}
We have:
$${\rm Out} (\pi_1 (M(P,c)))\cong {\rm Isom}(M(P,c))\cong \MCG^{\rm hom} (M(P,c)).$$ 
\end{cor}
\begin{proof}
Since $\Gamma$ is not elementary, the natural map
$\pi: \Isom (M(P,c))\to{\rm Out} (\pi_1 (M,c))$
is injective (in particular, there exists at most one isometry
in every homotopy class of self-homeomorphisms of $M(P,c)$).
By Theorem~\ref{group:isometry:teo} it follows that
$\pi$ is an isomorphism, so we are left to prove that
any self-homeomorphism of $M(P,c)$ is homotopic to an isometry.
 
Let $f$ be a self-homeomorphism of $M(P,c)$.
Then there exist an isomorphism
$f_\ast: \Gamma\to\Gamma$ and a $f_\ast$-equivariant
lift $\widetilde{f}:\matH^3
\to\matH^3$ of $f$ (\emph{i.e.}~ a lift of 
$f$ with $\widetilde{f} (\gamma (x))=
f_\ast (\gamma) (\widetilde{f} (x))$ for all
$\gamma\in\Gamma$, $x\in\matH^3$). By Theorem~\ref{group:isometry:teo},
an element
$g\in{\rm Isom} (\matH^3)$ exists such that $\varphi (\gamma)=
g\gamma g^{-1}$ for all $\gamma\in\Gamma$. We now define
$\widetilde{F} : \matH^3\times [0,1]\to\matH^3$ by setting 
$\widetilde{F}(x,t)=(1-t) \cdot \widetilde{f} (x)+
t\cdot g(x)$. Being $f_\ast$-equivariant, $\widetilde{F}$
projects onto a homotopy $F :M(P,c)\to M(P,c)$ between $f$
and an isometry, whence the conclusion.
\end{proof}

\begin{rem}\label{series:rem}
In~\cite{KMS}, a Kleinian group is by definition a discrete
torsion-free subgroup of $\Isom^+ (\matH^3)$ with non-empty 
\emph{discontinuity set} (see~\cite{Thu} for a definition). However, proofs
in~\cite{KMS} also work when dealing with Kleinian groups
with full-measure limit set~\cite{communication}.
Moreover, using Remark~\ref{exh:rem} one could easily prove
that if $(P,c)$ is an infinite decorated polyhedron, then every
finitely generated subgroup of $\pi_1 (M(P,c))$ has non-empty
discontinuity set, and is actually geometrically finite.
\end{rem}

\subsection{From homotopy to isotopy}
We are now left to prove that homotopic self-homeomorphisms of manifolds
arising from our construction are in fact isotopic.
We begin with the following:

\begin{defn}\label{endirr:defn}
Let $M$ be a non-compact manifold. We say that $M$ is \emph{end-reducible}
if there exist a compact set $W\subset M$ and a sequence
$\{\lambda_n\}_{n\in\matN}$ of simple loops in $M\setminus W$ with the 
following properties: any compact subset of $M$ intersects only a finite number
of $\lambda_i$'s, and each $\lambda_i$ is homotopically trivial
in $M$ and homotopically non-trivial in $M\setminus W$.
A non-compact manifold is \emph{end-irreducible} if it is not 
end-reducible. 
\end{defn}

\begin{rem}\label{def:rem}
The notion of  end-irreducibility was introduced
in~\cite{brown1}.
Our definition of end-irreducibility is proved to be
equivalent to the original one in~\cite[Lemma 3.1]{brown1}.
\end{rem}
 
It is shown in~\cite{brown2} that two homotopic self-homeomorphisms 
of an end-irre\-ducible
non-compact manifold are in fact isotopic.
Therefore, in order to conclude the proof of Theorem~\ref{main2:teo}
we only need to show the following:

\begin{prop}\label{homotopy2isotopy}
Let $(P,c)$ be an infinite decorated polyhedron.
Then $M(P,c)$ is end-irreducible.
\end{prop}
\begin{proof}
Suppose $M=M(P,c)$ is end-reducible, 
and let $W$ and $\{\lambda_n\}_{n\in\matN}$ 
be as in Definition~\ref{endirr:defn}. 
By Remark~\ref{exh:rem}, there exists a complete finite-volume
hyperbolic manifold with geodesic boundary $M_W\subset M$
containing $W$. Let
$C_1,\ldots,C_j$ be the toric cusps of $M$ meeting $M_W$. 
Up to rescaling, we can suppose that each $C_i$ is disjoint from $W$.
Let $T_i\subset M$ be the Euclidean torus bounding $C_i$, and set
\[
M'_W=(M_W\setminus \left(C_1\cup\ldots\ldots C_j)\right)\cup 
\left(T_1\cup\ldots\cup T_j\right).
\]
Since $M'_W$ is compact, 
there exists $N\gg 0$ such that
$\lambda_N\cap M'_W=\emptyset$. 
Suppose $\lambda_N\subset C_i$ for some $i$. 
Being non-trivial in $M\setminus W$, the loop $\lambda_N$
is non-trivial in $C_i$, whence in $M$, 
since the map $i_\ast: \pi_1 (C_1)\to \pi_1 (M)$
induced by the inclusion
is injective:
a contradiction.
We can thus suppose
$\lambda_N\subset M\setminus M_W$.
Since $\lambda_N$ is non-trivial in $M\setminus W$, it is non-trivial 
in $M\setminus M_W$. Moreover, 
each surface in $\partial M_W$ is geodesic, whence 
incompressible in $M$, and this 
readily implies that for any component
$L$ of $M\setminus M_W$ the map $i_\ast :\pi_1 (L)\to\pi_1 (M)$ 
induced by the inclusion
is injective. 
This gives $i_\ast (\lambda_N)\neq 1\in\pi_1 (M)$:
a contradiction.
\end{proof}

\bibliographystyle{amsalpha}
\bibliography{biblio}

\end{document}